\documentclass[12pt]{article}
\usepackage{palatino}
\usepackage{eucal}
\usepackage{amsfonts,amsmath,amssymb,latexsym}
\newtheorem{theorem}{Theorem}

\newtheorem{lemma}{Lemma}

\newtheorem{remark}{Remark}
\newcommand{\R}{\mathbb{R}}
\def\Ss{\mathbb{S}}

\def\Rr{\mathbb{R}}

\def\dd{\mathrm{d}}

\def\ae{\langle}
\def\ad{\rangle}

\def\<{\langle}
\def\>{\rangle}
\def\bea{\begin{eqnarray*}}
\def\eea{\end{eqnarray*}}
\def\be{\begin{equation}}
\def\ee{\end{equation}}
\newcommand{\po}{{\hspace*{-1ex}}{\bf .  }}
\def\proof{\noindent{\it Proof: }}
\def\qed{\ifhmode\unskip\nobreak\fi\ifmmode\ifinner\else
\hskip5 pt \fi\fi\hbox{\hskip5 pt \vrule width4 pt height6 pt
depth1.5 pt \hskip 1pt }}

\begin{document}

\title{Killing graphs with prescribed\\ mean curvature}
\author{M. Dajczer\thanks{Partially supported by Procad, CNPq and Faperj.},\,\, P. A. Hinojosa\thanks{Partially supported by
PADCT/CT-INFRA/CNPq/MCT Grant\,\#620120/2004-5.} \,and\, J. H. de
Lira\thanks {Partially supported by CNPq   and Funcap.}}

\date{}
\maketitle

\begin{abstract}
It is proved the existence and uniqueness of Killing graphs with
prescribed mean curvature in a large class of Riemannian manifolds.
\end{abstract}
\vspace{0.3cm}

{\small \noindent {\bf Keywords:} Killing graphs, prescribed mean
curvature, quasilinear elliptic PDE. \noindent {\bf MSC 2000:}
53C42, 53A10.}

\section{Introduction}

A basic strategy to obtain  hypersurfaces with prescribed mean
curvature in Euclidean space is to describe them non-parametrically
as solutions of a Dirichlet problem for a certain quasilinear
elliptic PDE. The solutions are then graphs over domains in totally
geodesic hypersurfaces of the ambient space. Classical references on
the subject are \cite{GT}, \cite{JS} and \cite{S}. Existence results
were obtained also for curved space forms in several formulations as
one may consult \cite{BS}, \cite{gs}, \cite{L}, \cite{LM}, \cite{NS}
and \cite{ni}. Recently, the cases of Riemannian and warped products
also deserved research efforts. For instance, see \cite{AD},
\cite{DR} and \cite{Sp}.

The mere possibility of applying similar analytical and geometrical
tools in all contexts seen above indicates that one must search for
general existence results in a large class of Riemannian spaces. The
presence of an isometric or conformal Killing vector field playing a
major role is certainly the distinguishing feature for all the
particular ambient spaces aforementioned.

In this paper, we deal with Riemannian spaces endowed with a Killing
vector field and define a notion of Killing graph in these general
ambients. We are then able to solve the corresponding Dirichlet
problem for prescribed mean curvature under hypothesis involving
domain data and the Ricci curvature. Since the ambient metric is
indeed warped in a sense we made precise later, the present article
should be considered as an extension and generalization of results
proved in \cite{DR}. In that sense, it is worth to mention that key
arguments in \cite{DR} do not longer work in the more general
framework considered here.

The Ricci curvature naturally arises in all {\it apriori} estimates
we made since they are based on comparison of geometric data. By its
turn, geometric comparison results follow from  the classical Jacobi
and Ricatti equations. The last one is used in Section~3 for
controlling the extrinsic geometry of  the barriers. 

We now explain more precisely the framework we are considering. Let
$M$ be a $(n+1)-$dimensional Riemannian manifold endowed with a
non-singular Killing vector field $Y$. We assume that the
distribution orthogonal to $Y$ is integrable. Then, the leaves are
easily seen to be  totally geodesic hypersurfaces. Let $\mathbb{P}$
be a fixed integral leaf and assume that the flow lines of the flux
$\Psi:\Rr\times\mathbb{P}\to M$ generated by $Y$ are complete. Given
a bounded domain $\Omega$ in $\mathbb{P}$, the Killing graph
$\Sigma$ associated to a function $u$ on $\bar\Omega$ is the
hypersurface
\[
\Sigma=\{q=\Psi(u(p),p): p\in\bar\Omega\}.
\]

Our results assure the existence of Killing graphs with prescribed
mean curvature $H$ and boundary data $\phi$. Here, the functions $H$
and $\phi$ are defined respectively on $\bar\Omega$ and $\Gamma$
where $\Gamma=\partial\Omega$. The problem of existence of such
Killing graphs is formulated in terms of a Dirichlet problem for a
divergence form elliptic PDE. The barriers we used for estimating
height and gradient of $u$ are the Killing cylinders. The Killing
cylinder $K$ over $\Gamma$ is ruled by the flow lines of $Y$ through
$\Gamma$. Therefore, we have
\[
K=\{q=\Psi(s,p):p\in\Gamma\}.
\]
In the sequel, the mean curvature of $K$ pointing inward is denoted
by $H_{\textrm{cyl}}$.\vspace{1ex}

We are now in condition to state the theorems proved in this paper.
We refer to Section 2 for the convention we used for the Ricci
tensor.

\begin{theorem}\po\label{main} Let $\Omega\subset \mathbb{P}$ be
a bounded domain with $C^{2,\alpha }$ boundary $\Gamma$. Suppose
that $H_{\rm{cyl}}\ge 0$ and that
$$
\inf_M {\rm Ric}_M \ge -n\inf_\Gamma H_{\rm{cyl}}^2.
$$
Let $H\in C^\alpha(\Omega)$ and $\phi\in C^{2,\alpha}(\Gamma)$ be
given. If $ \sup_\Omega|H|\le \inf_\Gamma H_{\rm{cyl}},$ then there
exists a unique function $u\in C^{2,\alpha}(\bar\Omega)$ satisfying
$u|_\Gamma=\phi$ whose Killing graph has mean curvature $H$.
\end{theorem}

If the Ricci tensor does not satisfy the assumption given in the
above result, then we may use comparison theorems with geodesic
spheres as barriers in order to prove the following result. In this
situation, we must impose certain condition either on the function
$H$ or in the size of the domain $\Omega$.

\begin{theorem}\po\label{main2} Let $\Omega\subset \mathbb{P}$ be
a bounded domain with $C^{2,\alpha }$ boundary $\Gamma$ contained
in~a normal geodesic disk with radius $r_0$. Suppose that $\inf_M
{\rm Ric}_M \ge -(n-1)k$ for some positive constant $k$. Let $H\in
C^\alpha(\Omega)$ and $\phi\in C^{2,\alpha}(\partial\Omega)$ be
given. \mbox{Assume} that $H_{\rm{cyl}}\ge 0$, that
$\sup_\Omega|H|\le \inf_\Gamma H_{\rm{cyl}}$ and that
$$
r_0 \le \frac{1}{\sqrt{k}}\coth^{-1}\frac{\sup_\Omega
|H|}{\sqrt{k}}.
$$
Then there exists a unique function $u\in C^{2,\alpha}(\bar\Omega)$
satisfying $u|_\Gamma=\phi$ whose Killing graph has mean curvature
$H$.
\end{theorem}

If the metric induced on $\mathbb{P}$ is rotationally symmetric then
we may take certain constant mean curvature spheres as barriers.
These spheres are rotationally invariant hypersurfaces and their
qualitative aspect is described by a flux formula. In what follows,
we denote $r=\textrm{dist}(p_0,\,\cdot\,)$ for some point
$p_0\in\mathbb{P}$.

\begin{theorem}\po\label{main3} Suppose that the induced metric in
$\mathbb{P}$ is of the form
\[
\dd r^2 +\xi^2(r)\dd\theta^2
\]
for a given function $\xi$, where $\dd \theta^2$ is the canonical
metric on the unit sphere $\Ss^{n-1}$. Let $\Omega\subset
\mathbb{P}$ be a bounded domain with $C^{2,\alpha }$ boundary
$\Gamma$ contained in a normal geodesic disk with radius $r_0$. Let
$H\in C^\alpha(\Omega)$ and $\phi\in C^{2,\alpha}(\partial\Omega)$
be given. \mbox{Assume}  that $H_{\rm{cyl}}\ge 0$,  that
$\sup_\Omega|H|\le \inf_\Gamma H_{\rm{cyl}}$ and that
$$
\sup_\Omega |H| \le -\frac{\varrho(r_0)\xi^{n-1}(r_0)}{\int_0^{r_0}
\varrho(r)\xi^{n-1}(r)\dd r}
$$
where $\varrho =|Y|^2$. Then there exists a unique function $u\in
C^{2,\alpha}(\bar\Omega)$ satisfying $u|_\Gamma=\phi$ whose Killing
graph has mean curvature $H$.
\end{theorem}

We emphasize that in the statements and proofs of the theorems we
may replace $M$ (after passing to the universal cover if necessary)
by the solid cylinder $\Psi(\mathbb{R}\times\bar\Omega)$ whose
boundary is $K$. Moreover, the hypothesis in Theorems 2 and 3 may be
understood either as restrictions on the size of the domains for
arbitrary mean curvature functions, or as upper bounds on $|H|$ in
the case of arbitrarily large domains. Finally, we remark that
standard regularity theorems imply that the results remain true for
continuous boundary data.

If the assumption $\sup_\Omega|H|\le \inf_\Gamma H_{\rm{cyl}}$ fails
at some point, we do not show here that the our results are no
longer true in the sense that there is a boundary data $\phi $ for
which no solution exists. Nevertheless, this is well known in some
cases, namely, the Euclidean space for standard graphs and has also
been proved for two types of Killing graphs in the hyperbolic space
in \cite{gs} and \cite{ni}.

This paper is organized as follows. In Section 2, we present the
basic geometric structure of the ambient spaces we are dealing with,
including calculations concerning the mean curvature of the Killing
cylinders and of the Killing graphs. In Section 3, we deduce the
height estimates. Sections 4 and 5 are devoted, respectively, to
boundary and interior gradient estimates. The last section presents
the proof of the theorems following the well-known continuity
method. An appendix contains the sketched proof of the flux formula.

We point out that ambient spaces with a Killing vector field
correspond in Lorentzian setting to the important notion of
stationary space-times on which the metric tensor is
time-independent. Prescribed mean curvature hypersurfaces work in
this context as Cauchy hypersurfaces for the initial value
formulation for the Cauchy problem in General Relativity. Thus, it
seems interesting to investigate the Lorentzian analogues to our
existence results.

The case of non-integrable orthogonal distributions corresponds to
Riemannian submersions other than the simple ones associated to
warped products. The particular situation of the three-dimensional
Heisenberg Lie group with a left-invariant metric was treated in
\cite{ALR}.  Two of the authors of the pressent paper address in
\cite{DL} the question for Riemannian submersions with
uni-dimensional totally geodesic fibers. The case of conformal
vector fields was studied by one of the authors in joint work with
F. Andrade in \cite{AL}.


\section{Killing graphs}

Let $M$ be a connected $(n+1)-$dimensional  Riemannian manifold
endowed with a non-singular Killing vector field $Y$. We denote the
metric and the Riemannian connection in $M$ by $\ae
\cdot\,,\cdot\ad$ and $\bar\nabla$, respectively. We assume that the
flow lines of $Y$ are complete and that the distribution
\[
p\in M \mapsto \{v\in T_p M : \ae v,Y\ad =0\}
\]
is integrable. Then, it is easy to verify that the integral leaves
are totally geodesic hypersurfaces.

Let $\mathbb{P}$ be such an integral leaf. The flux
$\Psi\colon\,\mathbb{R}\times\mathbb{P}\to M$ generated by $Y$ takes
isometrically  $\mathbb{P}$ to the leaves
$\mathbb{P}_s=\Psi_s(\mathbb{P})$, where $\Psi_s=\Psi(s,\,\cdot\,)$.
Given local coordinates $x^1,\ldots, x^n$ for $\mathbb{P}$, then
$s,x^1,\ldots, x^n$ are local coordinates for $M$ defined by
\[
q\in M \mapsto (s,x^1,\ldots, x^n) \quad\textrm{ if }\quad
q=\Psi(s,p),
\]
where $p\in \mathbb{P}$ is the point with coordinates $x^1,\ldots,
x^n$. The corresponding coordinate vector fields are
$$
\partial_s =\frac{\dd}{\dd s}\, \Psi(s,p) = Y(\Psi(s,p))
$$
and
$$
\partial_i(q)=\frac{\partial}{\partial x^{i}}\Psi(s,p)
=\Psi_{s_*}(p)\partial_i(p).
$$
The ambient metric in terms of these coordinates has components
$$
g_{00}=\<\partial_s,\partial_s\> = \varrho^2,\quad g_{0i}
=\<\partial_s,\partial_i\>=0
$$
and
$$
g_{ij}=\<\Psi_{s_*}\partial_i,\Psi_{s_*}\partial_j \>= \<\partial_i,
\partial_j \>= \sigma_{ij},
$$
where $\sigma_{ij}$ are the components of the metric in $\mathbb{P}$
in terms of the coordinates $x^i$. Observe that the components of
the metric do not depend on $s$. The gradient of the function $s$ is
$$
\bar{\nabla}s=g^{00}\,\partial_s = |Y|^{-2}\,Y =:f\,Y.
$$

Since the flow lines of $Y$  have constant geodesic curvature and
parallel curvature vector it is a standard fact that the solid
cylinder $\Psi(\mathbb{R}\times\bar\Omega)$ has a warped product
structure as
\[
\bar\Omega\times_\varrho\R
\]
whose metric is
\[
\sigma_{ij}\,\dd x^i \dd x^j + \varrho^2 \dd s^2.
\]
In fact, this is the  setting considered in \cite{DR}. Notice that
the relation between the Ricci curvatures of $M$ and $\mathbb{P}$
is determined by $\varrho=|Y|$.

Given a bounded $C^{2,\alpha}$ domain $\Omega$ on  $\mathbb{P}$ and
a smooth function $u$ on $\Omega$, we define the associated {\it
Killing graph} $\Sigma$ by
$$
\Sigma=\{q=\Psi(u(p),p):p\in \Omega\}.
$$
We may think of $\Sigma$ as the locus
$$
\Phi(s,p)=s-u(p)=0,
$$
where $u(p)=u(s,p)$. An orientation for $\Sigma$ is given at
$q=\Psi(u(p),p)$ by
\begin{eqnarray*}
\bar\nabla \Phi(q) \!\!&=&\!\! g^{00}\, \partial_s -g^{ij}(p)\,u_i\,
\partial_j(q)
=f\partial_s -\sigma^{ij}(p)\,u_i\, \Psi_{u(p)_*}\partial_j(p)\\
\!\!&=&\!\! f\,\partial_s - \Psi_{u(p)_*}\nabla u(p),
\end{eqnarray*}
where $\nabla$ denotes the Riemannian connection in $\mathbb{P}$ and
$$
\nabla u=\sigma^{ij}u_i\partial_j=u^j \partial_j
$$
is the gradient relatively to $\mathbb{P}$. Then \be\label{normal}
N= \frac{1}{W}\,\bar\nabla \Phi=\frac{1}{W}(f\,\partial_s -
\Psi_{_*}\nabla u) \ee defines a unit normal vector field along
$\Sigma$, where
$$
W^2 = f+|\nabla u|^2
$$
and $\Psi_*$ is a shorthand notation for $\Psi_{u(p)_*}$.

\subsection{Killing cylinder}

The {\it Killing cylinder\/} over $\Gamma=\partial\Omega$ is the
surface  ruled by the flow lines of $Y$ given by
$$
K=\{\Psi(s,p): s\in\mathbb{R},\; p\in \Gamma\}.
$$
If $s^1,\ldots, s^{n-1}$ are local coordinates for $\Gamma$, then
$s,s^1,\ldots, s^{n-1}$ are local coordinates for $K$. Let $\eta$ be
the unit inward normal vector along $\Gamma$ as a submanifold of
$\mathbb{P}$. We equally denote by $\eta$ the unit normal vector
field $\Psi_{s_*}\eta$ along $K$. Thus, we have
$$
\<\eta,\partial_s\>=0=\<\eta ,\partial_i\>.
$$
Since $\eta$ and $\partial_i$ are tangent to the totally geodesic
leaves $\mathbb{P}_s$, it results that
$$
\<\bar\nabla_{\partial_i}\partial_s, \eta\>=0.
$$
Hence $\partial_s$ is a principal  direction of $K$, and the
corresponding principal curvature is the geodesic curvature
$$
\kappa= f\<\bar\nabla_{\partial_s}\partial_s, \eta\>
$$
of the flow lines through $\Gamma$.

In the sequel, we deduce some useful properties of the distance
function $d=\textrm{dist}(\,\cdot\,,K)$ from $K$. We denote by
$\Gamma_\epsilon$ and $K_\epsilon$ the level sets $d=\epsilon$ in
$\mathbb{P}$ and $M$, respectively.  Thus, these level sets are
equidistant respectively from $\Gamma$ and $K$. It is immediate that
$K_\epsilon$ is the Killing cylinder over $\Gamma_\epsilon$. Since
$\Gamma$ is assumed to be $C^{2,\alpha}$, the function $d$ is also
$C^{2,\alpha}$ at points of $\Psi(\mathbb{R}\times\Omega_\epsilon)$,
where $\Omega_\epsilon\subset\Omega$ is a small tubular neighborhood
of $\Gamma$. Thus, we may define Fermi coordinates on
$\Psi(\mathbb{R}\times\Omega_\epsilon)$ as follows: for
$q\in\Psi(\mathbb{R}\times\Omega_\epsilon)$ we associate coordinates
$s^i,d$ by $q =\exp_{p}\, d\,\eta$ when $p=p(s,s^1,\ldots, s^{n-1})$
in $K$. Then \be\label{two} |\bar\nabla d|=1. \ee From this it
follows that \be\label{three} d^i d_{i;j}=0, \ee where
$d^i=g^{ij}d_j$ as usual. We also have
\[
\<\bar\nabla_{\partial_d}\bar\nabla d,
\partial_d\>= \frac{1}{2}\partial_d |\bar\nabla d|^2=0.
\]
Therefore,
\begin{eqnarray*}
\Delta d|_{d=\epsilon} \!\!& = &\!\!
\<\bar\nabla_{\partial_d}\bar\nabla d,
\partial_d\>+f\ae \bar\nabla_{\partial_s}\bar\nabla d,\partial_s\ad+
\sigma^{ij}\< \bar\nabla_{\partial_i}\bar\nabla d,\partial_j\>
\\
\!\!& = &\!\! f\ae \bar\nabla_{\partial_s}\bar\nabla
d,\partial_s\ad+ \sigma^{ij}\< \bar\nabla_{\partial_i}\bar\nabla
d,\partial_j\>.
\end{eqnarray*}
However $\bar\nabla d|_{d=\epsilon}=\partial_d$ is the unit inward
normal vector field $\eta_\epsilon$ to the equidistant cylinders
$K_\epsilon$. Denoting
\[
\kappa_\epsilon =f\ae \bar\nabla_{\partial_s}\partial_s,\bar\nabla
d\ad,
\]
we have at points of $K_\epsilon$ that \be \label{four} \Delta d
=-\kappa_\epsilon -\sigma^{ij}b_{ij}(\epsilon)
=-\kappa_\epsilon-(n-1)h_\epsilon, \ee where $b_{ij}(\epsilon)$ and
$h_\epsilon$ are the second fundamental form and the mean curvature,
respectively, of the hypersurface $\Gamma_\epsilon$ relatively to
the unit inward normal $\eta_\epsilon$. It follows that
$$
\Delta d |_{K_\epsilon}=-nH_{\textrm{cyl}}(\epsilon),
$$
where $H_{\textrm{cyl}}(\epsilon)$  is the mean curvature of
$K_\epsilon$ with respect to $\eta_\epsilon$.  Its Weingarten
operator is denoted by $A_\epsilon$. The mean curvature of $\Gamma$
and $K$ are denoted, respectively, by $h$ and $H_{\textrm{cyl}}$.

\begin{remark}\label{Omega_0}\po {\em All of the above calculations on
the distance function remain valid if we replace $\Omega_\epsilon$
by the larger  subset $\Omega_0$ in $\Omega$ consisting of the
points which can be joined to $\Gamma$ by a {\it unique} minimizing
geodesic. It was shown in \cite{LN} that in this set the function
$d$ has the same regularity as $\Gamma$. }\end{remark}

Throughout this paper, the ambient Ricci tensor in a given direction
$v$ is defined by
\[
\textrm{Ric}_M (v)=\sum_{i=1}^n\ae \bar R(e_i,v)v,e_i\ad,
\]
where $\bar R$ is the curvature tensor in $M$ and $e_1,\ldots,e_n,v$
is an orthonormal basis. We follow \cite{GT} or \cite{Sp} and use
the fact referred to in Remark \ref{Omega_0} on the distance
function in $\Omega_0$ for proving   in terms of the notation we
fixed above the following result.

\begin{lemma}\label{ricatti}\po
Assume that the  Ricci curvature satisfies ${\rm
Ric}_M\ge-n\inf_\Gamma H_{{\rm cyl}}^{\,2}$. Let $y_0\in\Gamma$ be
the closest point to a given point
$x_0\in\Gamma_\epsilon\subset\Omega_0$. Then, we have
$$
H_{{\rm cyl}}(\epsilon)|_{x_0}\ge H_{{\rm cyl}}|_{y_0}.
$$
\end{lemma}

\vspace{0.3cm}

\proof  We use local coordinates $s^0=s, s^1,\ldots, s^{n-1}$ as
defined above. At $d=\epsilon$ and since $\partial_d$ is the unit
speed of a geodesic, we have
\begin{eqnarray}
\label{ric1} -\frac{\dd}{\dd\epsilon} \<
A_\epsilon\partial_i,\partial_j\>
\!\!&=&\!\!\partial_d|_{d=\epsilon}
\<\bar\nabla_{\partial_i}\partial_d,\partial_j\>
=\<\bar\nabla_{\partial_d}\bar\nabla_{\partial_i}
\partial_d,\partial_j\> +
\<\bar\nabla_{\partial_i}\partial_d,\bar\nabla_{\partial_d}
\partial_j\>
\nonumber\\
\!\!& = &\!\! -\< \bar
R(\partial_i,\partial_d)\partial_d,\partial_j\> +\ae A_\epsilon
\partial_i,A_\epsilon \partial_j\ad
\end{eqnarray}
where $\bar R$ is the curvature tensor of $M$. On the other hand,
\begin{eqnarray}\label{ric2}\nonumber
\frac{\dd}{\dd \epsilon}\< A_\epsilon \partial_i,\partial_j\>\!\! &
= &\!\! \< \bar\nabla_{\partial_d}A_\epsilon \partial_i,
\partial_j\>+\ae
A_\epsilon\partial_i,\bar\nabla_{\partial_d}\partial_j\ad \\
\!\!&  = &\!\! \nonumber\<(\bar\nabla_{\partial_d}
A_\epsilon)\partial_i,\partial_j\> + \< A_\epsilon
\partial_j, \bar\nabla_{\partial_d}\partial_i\>
+\< A_\epsilon \partial_i,
\bar\nabla_{\partial_j}\partial_d\>\nonumber\\
\!\!& = & \!\!\<  A'_\epsilon \partial_i,\partial_j\> - 2\<
A_\epsilon
\partial_i, A_\epsilon \partial_j\>.
\end{eqnarray}
From (\ref{ric1}) and   (\ref{ric2}) we obtain the well-known
Ricatti equation
$$
A'_\epsilon - A_\epsilon^2 - \bar R_\epsilon=0
$$
where $\bar R_\epsilon=\ae
R(\cdot,\partial_d)\partial_d,\cdot\ad|_{d=\epsilon}$. Taking traces
we obtain
$$
n\frac{\dd}{\dd\epsilon}H_{\textrm{cyl}}(\epsilon) = \partial_d
\textrm{tr}A_\epsilon =\textrm{tr}\bar\nabla_{\partial_d} A_\epsilon
= \textrm{tr} \big(A_\epsilon^2 + \bar R_{\epsilon}\big) \ge
nH_{\textrm{cyl}}^2(\epsilon) + {\rm Ric}_M(\partial_d).
$$
 From our hypothesis on ${\rm Ric}_M$ we have that
$z=H_{\textrm{cyl}}(d)-\inf_\Gamma H_{{\rm cyl}}$ satisfies $$
z'(d)\geq z^2(d)-z^2(0)=(z(d)+z(0))(z(d)-z(0)).
$$
Thus $z'(d)\ge c(z(d)-z(0))$  in some interval $d\in [0,d_0>0]$ for
a constant $c>0$. It follows easily that
$H_{\textrm{cyl}}(\epsilon)$ does not decrease with increasing
$d$.\qed

\subsection{The mean curvature equation}

In what follows, we assume that the mean curvature $H$ of the
Killing graph $\Sigma$ is a function on $\bar\Omega$. Computing at
$q=\Psi(u(p),p)\in\Sigma$, we have
$$
nH(p) =-\textrm{tr}_{\Sigma} \bar\nabla N=-\textrm{div}_{\Sigma}N.
$$
Hence, if $e_1,\ldots,e_n$ is an orthonormal tangent frame at $q$ in
$\Sigma$, then
$$
-nH=\sum_i \<\bar\nabla_{e_i}N,e_i\>= \sum_i
\<\bar\nabla_{e_i}N,e_i\>+\<\bar\nabla_N N,N\>
=\textrm{div}_{M}\frac{\bar \nabla \Phi}{W}(q).
$$

Consider a normal coordinate frame $\partial_1,\ldots,\partial_n$
with $\sigma_{ij}=\delta_{ij}$ at $p\in\mathbb{P}$. Then, we have an
orthonormal frame $E_0=f^{1/2}\partial_s$,
$E_i(q)=\Psi_*\,\partial_i(p)$ at $q$. Using this frame and
(\ref{normal}), the divergence in the formula above becomes
$$
\textrm{div}_{M}\frac{\bar \nabla\Phi}{W}=
\<\bar\nabla\big(\frac{f}{W}\big),\partial_s\>
+\frac{f}{W}\textrm{div}_{M}\partial_s -f\<\bar\nabla_{\partial_s}
\frac{\Psi_*\nabla u}{W},\partial_s\>-
\<\bar\nabla_{E_i}\frac{\Psi_*\nabla u}{W},E_i\>.
$$
The Killing equation implies that $Y=\partial_s$ is divergence-free,
and that $f$ and $W$ do not depend on $s$. Thus,
\begin{eqnarray*}
\textrm{div}_{M}\frac{\bar \nabla\Phi}{W}\!\!&=&\!\!
-f\<\bar\nabla_{\partial_s} \frac{\Psi_*\nabla u}{W},\partial_s\>
- \sum_i\<\bar\nabla_{E_i}\frac{\Psi_*\nabla u}{W},E_i\>\\
\!\!&=&\!\!f\< \frac{\Psi_*\nabla u}{W},
\bar\nabla_{\partial_s}\partial_s\>-\sum_i\<
\bar\nabla_{\Psi_*\partial_i}\Psi_*\frac{\nabla u}{W},
\Psi_*\partial_i\>\\
\!\!&=&\!\!f\< \frac{\Psi_*\nabla u}{W},
\bar\nabla_{\partial_s}\partial_s\>-\sum_i\<\bar\nabla_{
\partial_i}\frac{
\nabla u}{W},\partial_i\>.
\end{eqnarray*}
Since $\Psi_*$ preserves the field $\partial_s$ and  $f$ is constant
along the flow lines, we have
$$
f(q)\< \frac{\Psi_*\nabla u}{W},
\bar\nabla_{\partial_s}\partial_s\>(q) =f(p)\<\frac{ \nabla u}{W},
\bar\nabla_{\partial_s}\partial_s\>(p).
$$
From these calculations it results that
$$
nH(p)= \sum_i\<\bar\nabla_{\partial_i}\frac{\nabla
u}{W},\partial_i\> -f\<\frac{ \nabla
u}{W},\bar\nabla_{\partial_s}\partial_s\>,
$$
where the expressions on both sides are now evaluated at $p\in
\mathbb{P}$. Since the $\partial_i$'s are orthonormal at $p$ and
$\mathbb{P}$ is totally geodesic, we may write \be \label{eqtn-1}
\textrm{div}_{\mathbb{P}}\big(\frac{\nabla u}{W}\big)- f\<\frac{
\nabla u}{W},\bar\nabla_{\partial_s}\partial_s\> -nH=0. \ee The
Killing equation implies that  the field
$\bar\nabla_{\partial_s}\partial_s$ is tangent to the leaf
$\mathbb{P}$.

On the other hand, it is easy to see that
$$
\<\bar\nabla f, \nabla u\>=
2f^2\<\bar\nabla_{\partial_s}\partial_s,\nabla u\>.
$$
Using this expression  and after some manipulation, one proves that
another way to write out (\ref{eqtn-1}) is
\begin{equation}
\label{eqtn-4}
\frac{1}{W}\big(\sigma^{ij}-\frac{u^iu^j}{W^2}\big)u_{i;j}
-\frac{1}{W^3}(f+W^2) \ae \nabla u,
\bar\nabla_{\partial_s}\partial_s\ad-nH=0
\end{equation}
where $u_{i;j}$  is the Hessian of $u$ in terms if the coordinates
$x^i$ in $\mathbb{P}$. Denoting
\begin{equation}
\label{a} a^{ij}(x,\nabla u)
=\frac{1}{W}\big(\sigma^{ij}-\frac{u^iu^j }{W^2}\big)
\end{equation}
and
\begin{equation}
\label{b} b(x,\nabla u)= -\frac{1}{W^3}(f+W^2)\<\nabla u,
f\bar\nabla_{\partial_s}\partial_s\>,
\end{equation}
the mean curvature equation becomes
$$
\label{Q} \mathcal{Q}[u] = a^{ij}u_{i;j}+  b-nH=0.
$$
The matrix $a^{ij}$ is positive-definite with eigenvalues
$$
\lambda=\frac{f}{W^3}\;\;\;\mbox{and}\;\;\; \Lambda=\frac{1}{W}
$$
with multiplicities $1$ and $n-1$ corresponding to the directions
parallel and orthogonal to $\nabla u$, respectively. Notice that
$\lambda\le\Lambda$ since $f\le W^2$ by definition.

Let $\phi$ be a $C^{2,\alpha}$ function on $\Gamma$. The Killing
graph of $\phi$ is a codimension two submanifold of $M$. Thus,
$\Sigma$ is a Killing graph with prescribed mean curvature $H$ and
prescribed boundary given by the graph of $\phi$ if and only if $u$
solves the Dirichlet problem
\begin{eqnarray}
\label{Dir} \mathcal{Q}[u]=0,\quad u|_\Gamma=\phi
\end{eqnarray}
for a quasilinear elliptic PDE. We may apply maximum and comparison
principles to (\ref{Dir}). Indeed, this follows from our hypothesis
that the function $H$ does not depend on $u$ (see, e.g., \cite{GT},
Ch. 10).

\section{Height estimates}

In this section, we obtain {\it apriori} $C^0$ estimates for
solutions of the Dirichlet problem (\ref{Dir}). We divided the
exposition in two cases concerning different  assumptions on the
ambient Ricci curvature.\vspace{1ex}

From now on, the distance function $d$ is regarded as the distance
from $\Gamma$ on the totally geodesic hypersurface $\mathbb{P}$.

\subsection{Killing cylinders as barriers}

We assume that the ambient Ricci curvature satisfies
\[
{\rm Ric}_M\ge-n\inf_\Gamma H_{{\rm cyl}}^{\,2},
\]
i.e., the hypothesis of  Lemma 1. In this case, we construct
barriers for $u$ in (\ref{Dir}) of the form \be\label{function}
\varphi(x)=\sup_\Gamma \phi + h(d(x)) \ee where the real function
$h$ will be chosen later. Along $\Omega_{0}$ we have \be\label{des}
\varphi_i = h' d_i \;\;\;\mbox{and}\;\;\; \varphi_{i;j}=h'' d_i d_j
+ h' d_{i;j}. \ee As in (\ref{two}) and (\ref{three}) above we have
$|\nabla d|^2 = d^i d_i =1$ and $d^id_{i;j}=0$. Moreover,
(\ref{four}) now reads as
$$
d^{i}_{\textrm{ };i}=\sigma^{ij}d_{i;j}=-(n-1)h_\epsilon.
$$
It is convenient to write $W^2 = f+ h'^2$. Then (\ref{eqtn-4})
yields \bea \mathcal{Q}[\varphi]+ nH  \!\!&=&\!\!
\frac{1}{W}\Big(\varphi^{i}_{\textrm{ };i} -\frac{\varphi^i
\varphi^j\varphi_{i;j}}{W^2}\Big) -\frac{1}{W^3}(f+W^2)\<
f\bar\nabla_{\partial_s}\partial_s, \nabla \varphi\>\\
\!\!&=&\!\! \frac{1}{W}\Big(-h'(n-1) h_\epsilon +h''-\frac{h'^2
h''}{W^2}\Big)-\frac{h'}{W^3}(f+W^2)\ae
f\bar\nabla_{\partial_s}\partial_s,\eta_\epsilon\ad\\
\!\!&=&\!\! \frac{f}{W^3}\big(h''-h'\ae
f\bar\nabla_{\partial_s}\partial_s,\eta_\epsilon\ad\big)-
\frac{h'}{W}\big((n-1)h_\epsilon +\ae
f\bar\nabla_{\partial_s}\partial_s,\eta_\epsilon\ad\big). \eea We
choose for (\ref{function}) the test function
$$
h=\frac{e^{CA}}{C}\big(1-e^{-Cd}\big)
$$
where $A>\textrm{diam}(\Omega)$ and $C$ is a positive constant to be
chosen later. Then,
$$
h'=e^{C(A-d)}\;\;\;\mbox{and}\;\;\; h''=-Ch'.
$$
Since the mean curvature of the equidistant cylinder $K_\epsilon$ is
given by
\[
nH_{\textrm{cyl}}(\epsilon) =(n-1)h_\epsilon+\kappa_\epsilon,
\]
we get
$$
\mathcal{Q}[\varphi]+nH =-\frac{fh'}{W^3}(C+\kappa_\epsilon)
-\frac{h'}{W}nH_{\textrm{cyl}}(\epsilon).
$$
Assuming $\sup_\Omega|H|\le \inf_\Gamma H_{\textrm{cyl}}$ and using
Lemma \ref{ricatti}, we obtain
$$
\mathcal{Q}[\varphi]+nH\le
-\frac{fh'}{W^3}(C+\kappa_\epsilon)-\frac{h'}{W}n|H|.
$$
Observe that $f/W^2\le 1$. Moreover, as $C\to \infty$ we have that
$$
\frac{h'}{W}=\frac{h'}{\sqrt{f+h'^2}}\to 1.
$$
Choosing $C\gg 0$ such that $C+\kappa_\epsilon>0$, we obtain
$$
\mathcal{Q}[\varphi]<-n(H+|H|)\le 0.
$$
We conclude that at points of $\Omega_0$ it holds that \bea
& & \mathcal{Q}[\varphi]<\mathcal{Q}[u]=0,  \\
& & \varphi|_\Gamma\ge u|_\Gamma. \eea

We now prove that $\varphi\ge u$ on $\bar\Omega$.  By contradiction,
assume that there exist points for which the continuous function
$\hat u:=u-\varphi$ satisfies $\hat u>0$. Hence $m:=\hat u(y)>0$ at
a maximum point $y\in \bar\Omega$ of $\hat u$. Choose a minimizing
geodesic $\gamma$ joining $y$ to $\Gamma$ for which the distance
$d=d(y,\Gamma)$ is attained. Thus, $\gamma(t)=\exp_{y_0}t\eta, \;
0\le t\le d$, starts from a point $y_0\in \Gamma$ with unit speed
$\eta$. Since $\gamma$ is minimizing, we have
$d(\gamma(t),\Gamma)=t$ and the function $\varphi$ restricted to
$\gamma$ is differentiable with $ \varphi'(\gamma(t))=e^{C(A-t)}$.
Since the maximum of $\hat u$ restricted to $\gamma$ occurs at
$t=d$, i.e., at the point $y$, one has that
$$
u'(\gamma(d))-\varphi'(\gamma(d))=\hat u'(\gamma(d))\ge 0.
$$
This implies that
$$
\<\nabla u(y),\gamma'(d)\>\ge \varphi'(\gamma(d))=e^{C(A-d)}>0.
$$
In particular $\nabla u(y)\neq 0$, and hence the level hypersurface
$$
S=\{x\in\Omega\cap B_r(y): u(x)=u(y)\}
$$
is regular for small radius $r$. Along $S$ we have
$$
\hat u (x)+\varphi(x)=\hat u (y)+\varphi(y)\ge \hat u
(x)+\varphi(y),
$$
and since $\varphi$ is an increasing function of $d(\,\cdot\,,
\Gamma)$ it follows that $d(x,\Gamma)\ge d(y,\Gamma)=d$. From this
we conclude that the points in $S$ are at a distance at least $d$
from $\Gamma$.  Since $S$ is $C^2$ it satisfies the interior sphere
condition: there exists a small ball $B_{\varepsilon}(z)$ touching
$S$ at $y$ contained in the side to which  $\nabla u(y)$ and
$\gamma'(d)$ points. Thus, the points of $B_{\varepsilon}(z)$
satisfy $u(x)\ge u(y)$, and hence
$$
\varphi(x)+m\ge u(x)\ge u(y)=\varphi(y)+m,\quad x\in B_\varepsilon
(z),
$$
where in the first inequality we used the definition of $m$. Again
because $\varphi$ is an increasing function of $d$, we have
$d(x,\Gamma)\ge d$ on $B_\varepsilon (z)$ and therefore this ball is
contained in the interior of $\Omega$ far away from $\Gamma$. This
allows us to extend the geodesic $\gamma$ through $B_\varepsilon
(z)$. We claim that the center $z$ of the ball is contained in this
extension. Otherwise, the broken line consisting of $\gamma$ and of
the radius in $B_\varepsilon(z)$ from $z$ to $y$ has length smaller
than {\it a} minimizing geodesic joining $z$ to $y_0\in\Gamma$ (for
a suitable small $\varepsilon$ such a geodesic must cross the level
hypersurface $S$ at a point $x\neq y$ at distance to $\Gamma$
greater than $d$). Thus, if there exists at least two distinct
minimizing geodesics joining $y$ to $\Gamma$, then the point $z$ is
contained in the extension of both geodesics after its intersection
at $y$. Choosing $\varepsilon$ sufficiently small, we see that this
configuration is not possible (the construction we made above
applies to both geodesics). This contradiction implies that the
maximum point $y$ belongs to $\Omega_0$. However, in this case,
$\hat u (y)\le 0$, a contradiction. We conclude that $u\le \varphi$
throughout $\bar\Omega$ and therefore $\varphi$ is a continuous
super-solution for the Dirichlet problem (\ref{Dir}).

In a similar way, we may construct lower barriers for $u$, that is,
continuous sub-solutions for (\ref{Dir}). It is clear that the
existence of these barriers implies the desired $C^0$ {\it apriori}
estimates.

\subsection{Geodesic spheres as barriers}

Next we assume only that the Ricci curvature has a finite lower
bound, that is,
\[
\textrm{Ric}_M \ge -(n-1)k,
\]
for some positive constant $k$.  In this case, we present a strategy
for obtaining height estimates for $u$. Our method  relies in a
Hessian comparison theorem and takes geodesic spheres as barriers.

We fix a point $p_0\in\mathbb{P}$ and consider the  function
$r=\textrm{dist}(p_0,\,\cdot\,)$.  Let ${\mathbb{H}^{n+1}(-k)}$ be
the hyperbolic space form with constant sectional curvature $-k$ and
$r_{\textrm{hyp}}$ a distance function on it. The Laplacian
Comparison Theorem (cf.\ \cite{SY}, p.\ 5, Corollary 1.1) yields
$$
\Delta r \le \Delta_{\mathbb{H}^{n+1}(-k)} r_{\textrm{hyp}}
$$
at corresponding (equidistant) points,  whenever $r$ is
differentiable. This implies that if we consider geodesic balls $B$
in $M$ (outside the cut locus of $M$) and $B_{\textrm{hyp}}$ in
$\mathbb{H}^{n+1}(-k)$ with same radius $r_0$, then the mean
curvatures of the respective geodesic spheres calculated with
respect to the gradient of the distances satisfy
$$
-H_{\partial B} \le -H_{\partial B_\textrm{hyp}}=-\sqrt{k}\coth
\sqrt{k}\,r_0.
$$
Thus, if we assume that our function $H$ satisfies
$$
r_0 \le \frac{1}{\sqrt{k}}\coth^{-1}\frac{\max_\Omega
|H|}{\sqrt{k}},
$$
then we have that
$$
|H|\le H_{\partial B}.
$$
Thus, if we suppose that the domain $\Omega$ is contained in a
geodesic disc of radius $r_0$ in $\mathbb{P}$, then the
corresponding geodesic sphere in $M$  is a barrier for $C^0$
estimates for the problem (\ref{Dir}). Indeed, it suffices to move
such a sphere along the flow lines of $Y$ and then apply the maximum
principle at a first tangency point between the graph $\Sigma$ and
the moving spheres.


\begin{remark}\po
{\em For constant mean curvature $H$ it was shown \cite{DR} that an
apriori height estimate exists if $\textrm{Ric}_M>-nH^2$.  This is
achieved by constructing a function that is subharmonic but only if
$H$ is constant.  Up to this estimate, all the aforementioned
results on Killing graphs in the Introduction for different ambient
spaces follow from Theorems \ref{main} and  \ref{main2} in this
paper. }\end{remark}

\subsection{CMC spheres as barriers}

Next we assume that the induced metric in $\mathbb{P}$ is
rotationally invariant. More precisely, the metric in $\mathbb{P}$
is of the form
$$
\dd r^2 + \xi^2(r) \dd \theta^2
$$
in terms of coordinates $(r,\theta)\in \R^+\times \Ss^{n-1}$, where
$\dd\theta^2$ denotes the usual metric in $\Ss^{n-1}$. We also
assume that $\varrho=\varrho(r)$, that is, the norm of the Killing
field does not depend on $\theta$. In this case, the ambient metric
is written in terms of cylindrical coordinates $s,r,\theta$ as
$$
\varrho^2(r) \dd s^2 +\dd r^2 + \xi^2(r) \dd \theta^2.
$$
and $M$ is a doubly-warped product with respect to warping functions
of the coordinate $r$.

A rotationally invariant hypersurface $\Sigma_0$ is parametrized by
an immersion  $\R\times \Ss^{n-1}\to M$ whose coordinate expression
is
$$
(u,\theta)\mapsto (s(u),r(u),\theta),
$$
where $u$ is the arc-length parameter of the profile curve
$\theta=\theta_0$. This means that $u$ is defined by $\varrho^2\dot
s^2 + \dot r^2 =1$ and that the induced metric in $\Sigma_0$ is
$$
\dd u^ 2 + \xi^2(u)\dd\theta^2.
$$

If $\Sigma_0$ has constant mean curvature $H_0$, then it satisfies a
first order equation given by the flux formula. In the flux formula
in the Appendix, we put $\Gamma$ as the geodesic circle of radius
$r=r(u)$ which is the intersection of $\Sigma_0$ and the leaf
$\mathbb{P}_s$, where $s=s(u)$. Thus we consider $D$ as the geodesic
disc in $\mathbb{P}_s$ with radius $r=r(u)$. The co-normal $\nu$ is
the unit velocity vector $\dot s
\partial_s +\dot r
\partial_r$. Finally, the Killing vector field $Y$ corresponds to
the coordinate vector field $\partial_s$. Thus, one has
\[
\< Y,\nu\> = \varrho^2(r(u))\dot s \,\,\textrm{ and }\,\, \< Y,N_D\>
=\< Y,\frac{Y}{|Y|}\>= \varrho(r(u)).
\]
Hence, plugging these expressions in the flux formula gives \bea
c\!\!&=&\!\! nH_0 \int_{D} \varrho+\int_{\Gamma} \dot s \varrho^2
=nH_0 \int_0^r \int_{\Ss^{n-1}} \varrho\xi^{n-1}\,\dd
r\dd\theta+\int_{\Ss^{n-1}} \dot s \varrho^2\xi^{n-1}\,\dd\theta.
\eea Since we  are integrating at a fixed value of $s$ and therefore
at a fixed value of $u$, we have
$$
nH_0\int_0^r \varrho\xi^{n-1}\dd r +\dot
s(r(u))\varrho^2(r)\xi^{n-1}(r) = \frac{c}{\omega_n}.
$$
The derivative of this expression with respect to $u$ gives a second
order ODE which characterizes CMC rotationally invariant
hypersurfaces in $M$. Compact solutions satisfy $r=0$ at the points
of maximum and minimum height, where $\dot s=0$. Hence, these
compact examples correspond to take $c=0$.

On the other hand, at maximum points for $r$ we have $\varrho^2 \dot
s^2=1$ and therefore
$$
nH_0\int_0^{r_0} \varrho(r)\xi^{n-1}(r)\dd r
+\varrho(r_0)\xi^{n-1}(r_0) = 0,
$$
where $r_0$ is the maximum value of $r$. Thus the mean curvature of
compact  rotational examples with maximum radius $r_0$ satisfy
$$
nH_0 =-\frac{\varrho(r_0)\xi^{n-1}(r_0)}{\int_0^{r_0}
\varrho(r)\xi^{n-1}(r)\dd r}=:-F(r_0).
$$
If we assume that $\Omega\subset \mathbb{P}$ is contained in a
geodesic disc with radius $r_0$, then these CMC spheres are barriers
for the height of a Killing graph with prescribed mean curvature
$H(x)$ satisfying
$$
n|H(x)|\le F(r_0).
$$
The value of $F(r_0)$ may be explicitly given in particular cases
such as $\mathbb H^{n+1}(-k)$ and $\mathbb{H}^n(-k)\times \R$.

\section{Boundary gradient estimates }

Our task now is to produce {\it apriori} gradient estimates for the
Dirichlet problem (\ref{Dir}). In order to do that, we use barriers
of the form $w+\phi$ along a tubular neighborhood $\Omega_\epsilon$
of $\Gamma$ as defined in Section 2.1. Here, $w=\psi(d(x))$ for some
real function $\psi$ to be chosen later and
$d=\textrm{dist}(\,\cdot\,,\Gamma)$. Moreover, the boundary data
$\phi$ was extended to $\Omega_\epsilon$ by $\phi(s^i,d)=\phi(s^i)$
for simplicity.

A simple estimate gives \bea \mathcal{Q}[w+\phi]\!\!&=&\!\!
a^{ij}(x,\nabla w+\nabla \phi)(w_{i;j}+\phi_{i;j})
+b(x,\nabla w+\nabla\phi)-nH\\
\!\!&\le&\!\! a^{ij}w_{i;j}+\Lambda|\phi|_{2,\alpha}+ b -nH, \eea
where $a^{ij}$ and $b$ are given by (\ref{a}) and (\ref{b}). Thus,
$$
a^{ij} w_{i;j}=\frac{1}{W}\Delta w -\frac{1}{W^3}
(w^i+\phi^i)(w^j+\phi^j)w_{i;j}.
$$
On one hand, we deduce from (\ref{two}) and (\ref{three}) the
expressions
\[
w^iw^jw_{i;j}=(\psi')^2 d^id^j(\psi'' d_i  d_j + \psi'd_{i;j})
=(\psi')^2\psi'' |\nabla d|^4=(\psi')^2\psi''
\]
and
\[
w^i\phi^jw_{i;j}=\psi' d^i \phi^j (\psi''d_i d_j + \psi'
d_{i;j})=\psi' \psi''|\nabla d|^2 d^j \phi_j =\psi'\psi''\<\nabla
d,\nabla \phi\>=0
\]
and
\[
\phi^i\phi^j w_{i;j}=\phi^i\phi^j(\psi'' d_i d_j+\psi'd_{i;j})
=\psi''\<\nabla d,\nabla
\phi\>^2-\psi'\phi^i\phi^jb_{ij}(\epsilon)=\!
-\psi'\phi^i\phi^jb_{ij}(\epsilon).
\]
In particular, we obtain
\[
\frac{1}{W}\psi''-\frac{1}{W^3}w^iw^jw_{i;j}
=\frac{\psi''}{W}\big(1-\frac{(\psi')^2}{W^2}\big)
=\frac{\psi''}{W^3}(f+|\nabla\phi|^2)
\]
since $W^2 = f + (\psi')^2 +|\nabla\phi|^2$ and
$$
(w^i\phi^j +w^j\phi^i +\phi^i\phi^j)w_{i;j}
=-\psi'\phi^i\phi^jb_{ij}(\epsilon).
$$
On the other hand, we deduce from (\ref{four}) that
\[
\Delta w = \psi'' +\psi'\Delta d = \psi'' -(n-1)\psi' h_\epsilon.
\]
We conclude that
\[
a^{ij}w_{i;j}= -\frac{\psi'}{W}(n-1)h_\epsilon +
\frac{\psi''}{W^3}(f+|\nabla\phi|^2) +\frac{\psi'}{W^3}
\phi^i\phi^jb_{ij}(\epsilon).
\]
A suitable expression for $b$ is \bea b \!\!&=&\!\!
-\frac{1}{W^3}(f+W^2) \<\nabla w+\nabla \phi,
f\bar\nabla_{\partial_s}\partial_s\>\\
\!\!&=&\!\! -\frac{\psi'}{W}\big(\frac{f}{W^2}+1\big)
\kappa_\epsilon -\frac{1}{W}\big(\frac{f}{W^2}+1\big)\ae
f\bar\nabla_{\partial_s}
\partial_s,\nabla \phi\ad
\eea since $\nabla w=\psi'\eta_\epsilon$ and
$\kappa_\epsilon=f\ae\bar\nabla_{\partial_s}\partial_s,\eta_\epsilon\ad$.
Therefore,
$$
\begin{array}{l}
\mathcal{Q}[w+\phi] \le -
\frac{\psi'}{W}\big((n-1)h_\epsilon+\kappa_\epsilon\big)
-nH+\Lambda|\phi|_{2,\alpha} -
\frac{\psi'f\kappa_\epsilon}{W^3}\vspace*{1.5ex}\\
\hspace*{.8in}
 - \frac{1}{W} \big(\frac{f}{W^2}+1\big)
\<f\bar\nabla_{\partial_s}\partial_s,\nabla
\phi\>+\frac{\psi''}{W^3}(f+|\nabla\phi|^2)
+\frac{\psi'}{W^3}\phi^i\phi^j b_{ij}(\epsilon).
\end{array}
$$
Finally, using that $\Lambda=1/W$, we get
$$
\begin{array}{l}
W^3\mathcal{Q}[w+\phi]\le -
n\psi'H_{\textrm{cyl}}(\epsilon)W^2-nHW^3
+|\phi|_{2,\alpha}W^2-\psi'f \kappa_\epsilon\vspace*{1.5ex}\\
\hspace*{6ex} - \ae f\bar\nabla_{\partial_s}\partial_s,\nabla
\phi\ad W^2 + \psi'\phi^i\phi^j b_{ij}(\epsilon)
-f^2\<\bar\nabla_{\partial_s}\partial_s,\nabla\phi\>
+\psi''(f+|\nabla\phi|^2 ). \vspace{1.5ex}
\end{array}
$$

Now define
\[
\psi (d) = \mu\ln (1+Kd)
\]
for certain positive constants $\mu$ and $K$ to be chosen later. We
have
$$
\psi'  = \frac{\mu K}{1+Kd}\quad\textrm{ and }\quad
\psi''=-\frac{1}{\mu}(\psi')^2.
$$
We choose $\mu$ in such a way that $\mu\to 0$ as $K\to \infty$. It
suffices to take
$$
\mu = \frac{C}{\ln (1+K)}
$$
for some positive constant $C$ to be chosen later. In this case, as
$K\to\infty$ one has
$$
\psi'(0)= \frac{CK}{\ln (1+K)}\to +\infty.
$$
It also holds that $\frac{\psi'}{W}\sim 1$ as $K\to \infty$. Thus,
at points of $\Gamma$ the last inequality (asymptotically)  becomes
\bea &  &  W^3\mathcal{Q}[w+\phi] \le -n(H_{\textrm{cyl}}+H)\psi'^3
-\frac{1}{\mu}(f+|\nabla\phi|^2)(\psi')^2\\
 &  & \hspace*{5ex}+\big(|\phi|_{2,\alpha}- \ae
f\bar\nabla_{\partial_s}\partial_s,\nabla\phi\ad\big) (\psi')^2-(f
\kappa +\phi^i\phi^j b_{ij})\psi'
-f^2\<\bar\nabla_{\partial_s}\partial_s,\nabla\phi\> \eea Therefore,
assuming that $H_{\textrm{cyl}}+H\ge 0$ and choosing $K$ large
enough, we assure that $\mathcal{Q}[w+\phi]<0$ on a small tubular
neighborhood $\Omega_\epsilon$ of $\Gamma$ and that $w+\phi\ge
u|_{\Gamma_\epsilon}+\phi$ on both boundary components. Therefore,
$w+\phi$ is a locally defined upper barrier for the Dirichlet
problem (\ref{Dir}). A lower barrier may be constructed in a similar
way.

\section{Interior gradient estimates}

The last step in providing {\it apriori} estimates for (\ref{Dir})
is to verify that $\nabla u$ satisfies a kind of maximum principle
for a third order equation obtained from differentiating
$\mathcal{Q}[u]=0$ and contracting the resulting equation with the
gradient itself.\vspace{1ex}

Using a suitable test function taken from \cite{TW} and Ricci
identities allows us to eliminate third derivatives and to obtain
global estimates for $|\nabla u|$ in terms of the height and
boundary $C^1$ estimates.

We suppose momentarily that $u\in C^3(\Omega)$. The usual regularity
theorems guarantee that the estimates we will obtain are also true
for a $C^{2,\alpha}$ function (see \cite{GT} and \cite{T}). Equation
(\ref{eqtn-4}) may be written as \be\label{eqtn-5}
A^{ij}\,u_{i;j}=B, \ee where \be\label{nine} A^{ij}(x,\nabla u) =
W^2\sigma^{ij}-u^i u^j \ee and \be\label{eleven} B(x,\nabla u)
=(f+W^2) \ae f\bar\nabla_{\partial_s}\partial_s,\nabla u\ad+nHW^3.
\ee Differentiating covariantly with respect to the metric
$\sigma_{ij}$ on $\mathbb{P}$ yields
$$
\nabla_k W^2 = \nabla_k f + \nabla_k (u^j u_j) =f_k+ u^j_{\textrm{
};k}u_j+ u^ju_{j;k}=f_k + 2u^j u_{j;k}.
$$
We obtain from (\ref{eqtn-5}) that
$$
\nabla_k B=((f_k + 2u^l u_{l;k}\big)\sigma^{ij}+ W^2
\nabla_k\sigma^{ij}-\nabla_k (u^i u^j))u_{i;j} +A^{ij}u_{i;jk}.
$$
Contracting with $u^k$ gives \be\label{der_eqn} u^k\nabla_k B=(f_k
u^k+ 2u^j u^k u_{j;k}) u^{i}_{\textrm{
};i}-2u^iu_{i;j}u^{j}_{\textrm{ };k} u^k +A^{ij}u^k u_{i;jk}. \ee
Following \cite{TW} we define the function
$$
\tau =e^{2Cu} v,
$$
where $v= u^i u_i$ is the squared norm of the gradient of $u$ and
$C>0$  a  constant to be chosen later. Then,
$$
\tau_i =e^{2Cu}(2Cvu_i + v_i).
$$
However,
$$
v_i=\nabla_i(u^j u_j) = u^j_{\textrm{ };i}u_j + u^j u_{j;i}=2 u^j
u_{j;i}.
$$
Thus, \be \label{one} \tau_i = 2e^{2Cu}(Cvu_i + u^j u_{j;i}). \ee If
$\tau$ achieves its maximum on $\Gamma$ then we have a bound for $v$
in $\Omega$ as desired. Hence, we may assume that the maximum is
attained at an interior point $p_0\in\Omega$ where we have
\be\label{eq_1} u^j u_{j;i}=-Cvu_i. \ee Differentiating (\ref{one})
yields \be\label{eight} \!\tau_{i;j}\!=\!2e^{2Cu}(2C^2u_i u_j v\!+
2C u_j u^k u_{k;i}+Cu_{i;j}v\! +2Cu_i u^{k}u_{k;j} \!+
u^{k}_{\textrm{ };j}u_{k;i}\!+u^k u_{k;ij}). \ee From the maximum
point criterion we have $\tau_{i;i}\le 0$, and from the ellipticity
of (\ref{eqtn-5}) that \be \label{max} A^{ij}\tau_{i;j}\le 0. \ee
From (\ref{eq_1}) we obtain at $p_0$ that $u_j u^j_{\textrm{};i}=
u^j u_{j;i}=-Cu_i v$. Therefore, \be \label{eq_3} u^j u_{j;i}u^i
=-Cv^2 \ee and \be\label{six} u^j u_{j;i}u^{i}_{\textrm{ };k}u^k
=C^2 v^2 u_k u^k = C^2v^3. \ee We may assume that $\nabla u$ is
non-singular in a neighborhood of $p_0$. Otherwise, we are done.
Hence, we choose local coordinates at $p_0$ asking $x^1$ to
parametrize the trajectories of $\nabla u$  and such that
$$
\partial_1 = \nabla u/|\nabla u|
$$
along the trajectory at $p_0$. The remaining coordinates parametrize
the level sets of $u$  and are assumed to be orthonormal and
geodesic at $p_0$.  Thus,  we have at $p_0$ that \be\label{ten}
u_1=|\nabla u|=v^{1/2}\;\;\;\mbox{and}\;\;\; u_j=0\;\;\mbox{for}\;\;
 j\neq 1.
\ee We also have $u_{j;i}=u^j_{{ };i}=u_{ji}$.  Then (\ref{eq_1})
gives \be\label{ten2} u_{1;1}=u_{11}=-Cv \;\;\;\mbox{and}\;\;\;
u_{1;i}=u_{1i}=0 \;\;\mbox{for}\;\;  i\neq 1. \ee If necessary, we
rotate the coordinates to assure that $u_{ij}$ is diagonal at $p_0$.
Then (\ref{eqtn-5}) at $p_0$ becomes \be\label{eq_4} (f+v)\Delta u =
(f+v)u_{ii} = B-vu_{11}= B-Cv^2. \ee We now use the Ricci identities
for the Hessian $u_{i;j}$ of $u$
$$
u_{i;jk}-u_{k;ij}= R_{ijkm}u^m,
$$
where $R$ is the curvature tensor in $\mathbb{P}$. Thus, we obtain
that \be\label{aes}
\begin{array}{l}
A^{ij}u^ku_{i;jk}=
((f+v)\sigma^{ij}-u^i u^j)u^k u_{i;jk}\vspace{1.5ex}\\
\hspace{.7in} =((f+v)\sigma^{ij}-u^i u^j)u^k(u_{k;ij}
+ R_{ijkm}u^m)\vspace{1.5ex}\\
\hspace{.7in} =A^{ij}u^ku_{k;ij}+(f+v)\sigma^{ij}R_{ijkm}u^k u^m
\end{array}
\ee since $R_{ijkm}u^iu^ju^ku^m=0$. Moreover, we have from
(\ref{eq_3}) and (\ref{eq_4}) that \be\label{seven} (f_ku^k +2u^j
u^k u_{j;k})u^i_{{ };i}=(f_k u^k-2Cv^2)D_1 \ee where
$D_1=:(B-Cv^2)/(f+v)$.  From (\ref{der_eqn}), (\ref{six}),
(\ref{aes}) and (\ref{seven}) we deduce \be u^k \nabla_k B = (f_k
u^k -2Cv^2)D_1 -2C^2 v^3+A^{ij}u^ku_{k;ij} +(f+v)\sigma^{ij}
R_{ijkm}u^ku^m. \ee We conclude that \be\label{third_der_1}
\!A^{ij}u^ku_{k;ij} =\!-\!f_k u^kD_1-(f+v)\sigma^{ij}\bar
R_{ijkm}u^k u^m +2C^2v^2f +2Cv^2D_2 +u^k\nabla_kB \ee where we
denote $D_2=(B-Cf^2)/(f+v)$.

Next we write down the maximality condition (\ref{max}) for $\tau$
at the point $p_0$.  A straightforward computation using
(\ref{nine}), (\ref{ten}) and (\ref{ten2}) in (\ref{eight}) yields
\bea 0\!\!&\ge&\!\! A^{ij}(2C^2u_iu_jv + 2Cu_ju^ku_{k;i} + Cu_{i;j}v
+ 2Cu_i u^k u_{k;j} +u^{k}_{\textrm{};j}u_{k;i}
+ u^k u_{k;ij})\\
\!\!&=&\!\! (f+v)(Cvu_{ii}+u_{jj}u_{jj})-2C^2 v^2
f+A^{ij}v^{1/2}u_{1;ij}. \eea Using that $u_{jj}u_{jj}\ge C^2v^2 $
and (\ref{eq_4}) we obtain \be\label{third_der_2}
A^{ij}v^{1/2}u_{1;ij}\le C^2v^2f -CBv. \ee Combining the expressions
(\ref{third_der_1}) and (\ref{third_der_2}) we get
$$
-f_1v^{1/2}D_1 +(f+v) \sigma^{ij} R_{ikjm}u^ku^m+C^2v^2f
+2Cv^2D_2+CBv+v^{1/2}\nabla_1 B\le 0.
$$
Let $R$ be the minimum eigenvalue of the Ricci tensor in
$\mathbb{P}$ in the direction of $\nabla u$. Thus, we have $Rv\le
\sigma^{ij} R_{ikjm}u^ku^m$ and  multiplying the above expression by
$f+v$ yields \be\label{ineq-1} \!Rv(f+v)^2+C^2 v^2 f(v-f)
+CBv(f+3v)+f_1v^{1/2}(Cv^2 -B) +(f+v)v^{1/2}\nabla_1 B\le 0. \ee
Next we compute the last term on the left hand side of
(\ref{ineq-1}). From (\ref{eleven}) we obtain \be\label{gamma}
B(x,\nabla u,v) =(2f^2 + fv)\Gamma_{00j}u^j + nH(f+v)^{3/2} \ee
where $\Gamma_{00j}=\sigma_{ij}\Gamma^i_{00}$. Differentiating
$B=B(x,\nabla u,v)$ we have
$$
\nabla_k B  =B_k + B_{u^j}u^j_{\textrm{ };k}+B_v v_k.
$$
We obtain at $p_0$ that
$$
v^{1/2}\nabla_1 B= B_1v^{1/2}-CB_{u^1}v^{3/2}-2B_vCv^2.
$$
Replacing this in (\ref{ineq-1}) yields \bea\label{ineq-2} & &
\nonumber Rv(f+v)^2 + C^2 v^2f(v-f) +(f+v)B_1
v^{1/2}+f_1v^{1/2}(Cv^2 -B)
\\
& &  \hspace*{8ex}+ Cv((f+3v)B-(f+v)(B_{u^1}v^{1/2}+2B_vv))\le 0.
\eea Using (\ref{gamma}) we have at $p_0$ that
$$
B_1=(4f+v)f_1\Gamma v^{1/2} +(2f^2+fv)\nabla_1\Gamma v^{1/2} +
nH_1(f+v)^{3/2}+\frac{3}{2}nH(f+v)^{1/2}f_1
$$
and
$$
B_{u^1}=(2f^2+fv)\Gamma
$$
and
$$
B_v = f\Gamma v^{1/2}+\frac{3}{2}nH(f+v)^{1/2},
$$
where $\Gamma=:\Gamma_{001}$. In particular,
$$
\!B_1 v^{1/2}\!=\!2f(2f_1\Gamma \!+f\nabla_1\Gamma)v +(f_1\Gamma
+f\nabla_1\Gamma)v^{2} +
n(f+v)^{1/2}v^{1/2}(H_1(f+v)\!+\frac{3}{2}Hf_1)
$$
and
$$
(f+3v)B-(f+v)(B_{u^1}v^{1/2}+2B_vv) = 2f^2\Gamma
v^{3/2}+nHf(f+v)^{3/2}.
$$
We may assume that $\max_\Omega f\le v$ since, otherwise, we are
done. Then, we have  $f+v\le 2v$ and \bea B_1 v^{1/2}
\!\!\!&\ge&\!\!\! (4ff_1\Gamma +2f^2\nabla_1\Gamma
\!-\!(3/\sqrt{2})n|Hf_1|)v\!+\! (f_1\Gamma +f\nabla_1\Gamma
-2^{3/2}n|H_1|)v^{2}\\
\!\!\!&=&\!\!\! A_1 v + A_2 v^2. \eea Moreover, we have
$$
(f+3v)B-(f+v)(B_{u^1}v^{1/2}+2B_v v)\ge (2f^2\Gamma
-2^{3/2}n|H|f)v^{3/2}=A_3 v^{3/2}
$$
and \bea f_1 v^{1/2}(Cv^2 -B) \!\!&=&\!\!
f_1v^{1/2}(Cv^2-(2f^2+fv)\Gamma v^{1/2}
-nH(f+v)^{3/2})\\
&\ge& -f_1(2f^2+fv)\Gamma v
-2^{3/2}n|Hf_1|v^2+Cf_1v^{5/2}\\
\!\!&=&\!\!C_1 v + C_2 v^2 + C f_1 v^{5/2}. \eea Replacing in
(\ref{ineq-2}) and grouping equal powers of $v$, we obtain
$$\begin{array}{l}
(C^2f+ R+A_2)v^3 + C(A_3+f_1)v^{5/2}+ (2 R f-C^2f^2 +
A_1+fA_2+C_2)v^{2}\\
\hspace*{10ex} +( Rf^2+A_1 f + C_1)v \le 0.
\end{array}
$$
Notice that all constants in the expression above depend only on $H$
and $\Omega$.  We choose $C$ large enough so that the coefficient of
$v^3$ in the left hand side is positive. With this choice we
conclude that $v\le \tilde C$ at $p_0$ for some constant $\tilde C$
depending on the height estimates for $u$ and on the data $H$
and~$\Omega$.


\section{Proof of the theorems}

We apply the well-known continuity method to the family of Dirichlet
problems
\begin{eqnarray*}
\mathcal{Q}_\sigma [u]=0,\quad u|_\Gamma = \sigma\phi,
\end{eqnarray*}
where $\sigma\in[0,1]$ and
$$
\mathcal{Q}_\sigma[u]=a^{ij}u_{i;j} + b -n\sigma H.
$$
The subset $I$ of $[0,1]$ consisting of values of $\sigma$ for which
the above Dirichlet problem has a $C^{2,\alpha}$ solution is
non-empty, since $0\in I$. The hypothesis $H_u=0$ implies that $I$
is open. This follows from a standard application of the implicit
function theorem. The closedness of $I$ follows from the {\it
apriori} estimates we had proved. Thus, the continuity method
assures that $1\in I$.

In order to prove uniqueness, it suffices to reproduce the proof
presented in  \cite{DR}.\qed\vspace{1ex}

We point out that our existence results still hold if $\phi$ is only
assumed continuous. We may approximate $\phi$ uniformly by smooth
boundary data and use the interior gradient estimate to obtain
strong convergence on compact subsets of $\Omega$. A local barrier
argument shows that the limiting solutions achieves the given
boundary data.

\section{Appendix: a flux formula}

Consider a hypersurface $\Sigma$ in $M$ with boundary $\Gamma$ and
another hypersurface $D$ such that $M\cup D$ bounds a domain $U$ in
$M$. Let $H$ be a function that is constant along the flow lines of
$Y$. Choose an orientation on $M\cup D$ given by unit normal vector
fields $N$ along $\Sigma$ and $N_D$ along $D$. By the Killing
equation, we have
$$
\textrm{div}_M HY = H\textrm{div}_M Y + \<\bar\nabla H,Y\>=0.
$$
Thus, applying Stokes theorem to the domain $U$, one has
$$
0 = \int_U \textrm{div}_M HY= \int_\Sigma H\ae Y,N\ad +\int_D H\ae
Y,N_D\ad.
$$
However, if we consider an orthonormal adapted frame
$N,e_1,\ldots,e_n$ to $\Sigma$ we have from the Killing equation
$$
0=\sum_i \<\bar\nabla_{e_i} Y, e_i\> = \sum_i \<\bar\nabla_{e_i} Y^T
, e_i\>+ \sum_i \ae\bar\nabla_{e_i}N, e_i\> =\textrm{div}_\Sigma Y^T
-nH \ae Y,N\ad.
$$
Thus, applying Stokes theorem to $\Sigma$ we obtain
$$
n\int_\Sigma H\ae Y,N\ad=\int_\Gamma \ae Y,\nu\ad
$$
where $\nu$ is the exterior unit co-normal of $\Sigma$ along
$\Gamma$. The two expressions we obtained above may be gathered at
the {\it flux formula}
$$
n\int_D H\ae Y,N_D\ad + \int_\Gamma\ae Y,\nu\ad=0.
$$
An useful variant of this reasoning is to consider a hypersurface
$\Gamma$ in $\Sigma$ not homologous to zero in $\Sigma$ such that
$\Gamma=\partial D$. Then we have that
$$
\label{flux} n\int_DH\ae Y,N_D\ad + \int_\Gamma\ae Y,\nu\ad=c,
$$
where $c$ is a constant depending only on the homology class of
$\Gamma$.

{\renewcommand{\baselinestretch}{1} \hspace*{-20ex}\begin{tabbing}
\indent \= Marcos Dajczer\\
\> IMPA \\
\> Estrada Dona Castorina, 110\\
\> 22460-320 -- Rio de Janeiro -- Brazil\\
\> marcos@impa.br\\
\end{tabbing}}

\vspace*{-5ex}

{\renewcommand{\baselinestretch}{1} \hspace*{-20ex}\begin{tabbing}
\indent \=Pedro A. Hinojosa \\
\> Departamento de Matematica \\
\> Universidade Federal da Paraiba - CCEN - DM \\
\> Cidade Universitaria\\
\> 58051-900 -- Jo\~ao Pessoa -- PB -- Brazil\\
\> hinojosa@mat.ufpb.br
\end{tabbing}}

\vspace*{-2ex}

{\renewcommand{\baselinestretch}{1} \hspace*{-20ex}\begin{tabbing}
\indent \= Jorge Herbert de Lira\\
\> Departamento de Matematica, \\
\> Universidade Federal do Ceara,\\
\> Bloco 914 -- Campus do Pici\\
\> 60455-760 -- Fortaleza -- Ceara -- Brazil\\
\> jherbert@mat.ufc.br
\end{tabbing}}

\end{document}